\documentclass[11pt,letterpaper]{amsart}
\usepackage{amssymb, amscd, amsfonts, euler, epic, xy,charter}
\xyoption{all}

\def\CC{{\mathbb C}}  
 
\def\FF{{\mathbb F}} 
 
\def\HH{{\mathbb H}}

\def\PP{{\mathbb P}}
\def\QQ{{\mathbb Q}} 
\def\RR{{\mathbb R}}

\def\ZZ{{\mathbb Z}}

\def\hor{{\rm hor}}

\def\ver{{\rm vert}}

\def\Acal{{\mathcal A}}
\def\Bcal{{\mathcal B}}
\def\Ccal{{\mathcal C}}
\def\Dcal{{\mathcal D}}
\def\Ecal{{\mathcal E}} 
\def\Fcal{{\mathcal F}} 
\def\Hcal{{\mathcal H}}

\def\Kcal{{\mathcal K}}
\def\Lcal{{\mathcal L}}

\def\Ocal{{\mathcal O}}
\def\Pcal{{\mathcal P}}

\def\Scal{{\mathcal S}}

\def\la{\langle}
\def\ra{\rangle}

\def\mfrak{\mathfrak{m}}

\def\Ubar{{\overline{U}}}

\def\ss{{\rm ss}}

\def\End{\operatorname{End}}
\def\pt{{\scriptscriptstyle\bullet}}

\def\half{{\tfrac{1}{2}}}

\newcommand\ad{\operatorname{ad}}

\newcommand\Hom{\operatorname{Hom}}

\newcommand\ord{\operatorname{ord}}

\newcommand\res{\operatorname{Res}}

\newcommand\spfrak{\operatorname{\mathfrak{sp}}}

\newcommand\sym{\operatorname{Sym}}

\newcommand\trace{\operatorname{trace}}

\newcommand\Sp{\operatorname{Sp}}
\newcommand\Mp{\operatorname{Mp}}

\newcommand\sign{\operatorname{sign}}

\newtheorem{theorem}{Theorem}[section]
\newtheorem{lemma}[theorem]{Lemma}
\newtheorem{proposition}[theorem]{Proposition}

\newtheorem{corollary}[theorem]{Corollary}

\theoremstyle{definition}
\newtheorem{definition}[theorem]{Definition}

\newtheorem{example}[theorem]{Example}

\theoremstyle{remark} 
\newtheorem{remark}[theorem]{Remark}

\author{Alex Boer }
\author{Eduard Looijenga}
\email{boer@math.uu.nl, looijeng@math.uu.nl}
\address{Mathematisch Instituut\\
Universiteit Utrecht\\
P.O.~Box 80.010, NL-3508 TA Utrecht\\
Nederland}

\title{On the unitary nature of abelian conformal blocks}
\begin{document}
\subjclass[2000]{81R10, 17B65}
\keywords{Fock representation, WZW connection, abelian conformal block}
\begin{abstract} 
We determine  the projectively flat unitary structure on abelian conformal blocks
in terms of WZW-data.
\end{abstract}

\maketitle

\section*{Introduction}
Although it is well-known that an  abelian conformal block is the same thing as a complete space of theta functions,  the equivalence of these two notions has, we believe,  not yet been sufficiently explicated. This paper is a contribution towards clearifying that 
equivalence, where the emphasis is on the unitary structure. Our motivation is that the classical side of the equation, the theory of theta functions, furnishes a canonical inner product (which originates from a polarized Hodge structure), whereas at the other side no-one yet succeeded in writing down an inner product on a general (nonabelian) conformal block that makes the WZW-connection unitary. Since a generic nonabelian conformal block on a given Riemann surface embeds in an abelian one on (what is called) a \emph{cameral} cover of that Riemann surface, our goal was, by way of exercise, to  transfer the inner product from one side to the other in the hope that the result would suggest what the inner product in the nonabelian case would  look like. We did that exercise  and we believe that the outcome, namely the inner product in terms of  WZW data, is sufficiently interesting to justify publication. For example, the inner product that we find turns out \emph{not} to be one that comes from the Fock representation. This implies that  in the nonabelian case we cannot expect the inner product to comes from one on an integrable  highest weight representation of the associated affine Lie algebra. But our work also suggests that the Hodge theory of covers of a Riemann surface will enter (perhaps in the guise of iterated integrals). We wish to emphasize that our treatment  is entirely algebraic (and hence, we are tempted to add, 
simple) to the extent that at no point we encounter (let alone write down) a theta function.
Moreover, the discussion is largely self-contained.

Abelian conformal blocks are also the topic of a recent paper by Andersen and Ueno
\cite{andersenueno}. 

Here is a brief description of contents by section. The first three  are devoted
to the theory of Fock representation in relation to  weight one Hodge structures.
Section one reviews the construction of the Fock representation starting from a symplectic vector space. In the next section we explain how a weight one Hodge structure polarized by the symplectic form defines  an inner product on its associated Fock representation. In the third section we do this with parameters and we identify the unitary connection of the Fock bundle attached to variation of polarized Hodge structure. We do the basic curvature computation: we check that the curvature of the Fock bundle is scalar and half that of the Hodge determinant bundle.

In section four we begin making the link with the WZW-connection. We first review the definition (in the spirit of \cite{looijeng:cfr}) and introduce the notion \emph{subalgebra of Fock type} (the ambient algebra is a local field or a generalization thereof). We found this notion,  which is also implicit in \cite{looijeng:cfr},  quite useful,  because it
isolates what is needed to make the transition from the WZW-connection to the
Fock connection. The main results are Propositions \ref{prop:fockiso} and \ref{prop:A}.
The last two sections have geometric content: we apply the preceding to a punctured Riemann surface resp.\ to  a family thereof. Among other things, we  relate the Fock connection with the residue pairing to the Gauss-Manin connection with the intersection pairing.
The main result here is Theorem \ref{thm:fockiso}. 
\\

It is a pleasure to acknowledge helpful conversations we had  with Hessel Posthuma during the preparation of this paper.
\\
\tableofcontents
\textbf{Convention.}
Suppose $R$ is a commutative ring with unit and $V$ a $R$-module. The symmetric group $\Scal_k$ acts on the $k$-fold tensor power
$\otimes_R^kV$ and we denote by $\sym^R_k(V):=(\otimes_R^k V)_{\Scal_k}$ the quotient $R$-module of $\Scal_k$-covariants and by $\sym_R^k(V):=(\otimes_R^k V)^{\Scal_k}$ the submodule of $\Scal_k$-invariants. Notice that $\sym^R_\pt(V)$ is naturally a commutative $R$-algebra (likewise, $\sym_R^\pt(V)$ is  a cocommutative
$R$-coalgebra). If $R\supset\QQ$ and $V$ is a free $R$-module, then the natural map $\sym_R^k(V)\to\sym^R_k(V)$ is an isomorphism of $R$-modules. 

\section{Fock representation for a symplectic vector space}
Let $k$ be a field of characteristic zero and let $H$ be a  finite dimensional $k$-vector space endowed with a nondegenerate antisymmetric form $(\; ,\: )$. This makes $H$ a  symplectic vector space. We note that under the identification 
\[
E: H\otimes H\to \End (H),\quad E(a\otimes b)(x)=2(b,x)a
\]
the subspace of symmetric tensors $\sym^2(H)$ gets identified with the Lie algebra $\spfrak (H)$ of the Lie group $\Sp (H)$.

The Heisenberg (Lie) algebra $\hat{H}$ associated to $H$ has $H\oplus k\hbar$ as underlying vector space and bracket 
$[a+\lambda \hbar, b+\mu\hbar]=(a,b)\hbar$. Its universal enveloping algebra
$U\hat H$ is the quotient of the tensor algebra of $\hat H$ by the 2-sided ideal generated by the elements $a\otimes b- b\otimes a-(a,b)\hbar$. If we stipulate that $H$ has degree $1$ and $\deg\hbar =0$, then this ideal consists of even degree tensors and so $U\hat H$ has a natural $\ZZ/2$-grading as a $k[\hbar]$-algebra:
\[
U\hat H=(U\hat H)_0\oplus (U\hat H)_1.
\]
It is clear that the image of $H\otimes H$ in $U\hat H$ maps under the projection 
$U\hat H\to UH=\sym_\pt(H)$ onto $\sym_2(H)$ with kernel $\CC\hbar$. We therefore 
denote this image by $\widehat{\sym}_2(H)$. Accordingly, the image of $\alpha\in H\otimes H$  in $\widehat{\sym}_2(H)$ is denoted by $\hat \alpha$.  In particular, 
$\widehat{E^{-1}}$ takes values in $(U\hat H)_0$.

\begin{lemma}\label{lemma:bracket}
The map $\tau: \spfrak (H)\to U\hat H[\hbar^{-1}]$,  $\tau (A):=\hbar^{-1}\widehat{E^{-1}(A)}$, is a 
Lie homomorphism: $[\tau(A) ,\tau(B)]=\tau([A,B])$ and $\ad (\tau (A))$ acts on 
$U\hat H[\hbar^{-1}]$ as $A$:
\[
[\tau (A),x_1\circ x_2\circ\cdots \circ x_k]=
\sum_{i=1}^k x_1\circ x_2\circ\cdots\circ A(x_i)\circ\cdots \circ x_k.
\]
\end{lemma}
\begin{proof}
By linearity it suffices to check this in case $A$ resp.\ $B$ is of the form
$E(a\otimes a)$ resp.\ $E(b\otimes b)$. Then
we have for any $x\in H$ in $U\hat H$ the following identity
\[
(a\circ a)\circ x= x\circ (a\circ a) +2(a,x)a\circ\hbar =
x\circ (a\circ a) + A(x)\circ\hbar
\]
which yields the last part of the lemma with induction. It also follows that
\begin{align*}
(a\circ a)\circ (b\circ b)&=(a\circ a\circ b)\circ b=\big(b\circ (a\circ a) +2(a,b)a\big)\circ \hbar)\circ b\\
&= b\circ (a\circ a\circ b) +2(a, b)a\circ b\circ\hbar\\
&= b\circ\big( b\circ (a\circ a) +2(a,b)a\circ \hbar\big)+2(a, b) a\circ b\circ\hbar\\
&=( b\circ b)\circ (a\circ a) +2(a,b) \big(b\circ a+a\circ  b\big)\circ\hbar.
\end{align*}
On the other hand, if $x\in H$, then
\begin{multline*}
[A,B](x)=4( a,b)(b, x)a-4( b,a)(a, x)b\\
=4( a,b)\big( (b, x)a+(a,x)b\big)=2(a,b)E(a\otimes b+b\otimes a)
\end{multline*}
and  so the lemma follows.
\end{proof}

Let $F$ be a maximal totally isotropic subspace $F\subset H$ relative to $(\; ,\: )$. Then 
$\hat F:=F\oplus k\hbar$ is an abelian subalgebra and hence
the linear form $\hat F\to k$ given by the coefficient of $\hbar$ may be regarded as a character of $\hat F$. The \emph{Fock representation of $\hat H$ defined by $F$} is the 
representation of $\hat H$ induced by this character:
\[
\FF (H,F):= U\hat H\otimes_{U\hat F}\CC.
\]
It is obtained from $U\hat H$ by dividing out by the right ideal generated by 
$F$ and $\hbar-1$. That ideal is $\ZZ/2$-graded and  so $\FF (H,F)$  inherits a $\ZZ/2$-grading from $ U\hat H$: $\FF (H,F)=\FF (H,F)_0\oplus \FF (H,F)_1$.
We shall denote by $v_o\in \FF (H,F)_0$ the generator (the image of $1\otimes 1$) and by $\rho: U\hat H\to \End(\FF (H,F))$ the action of $U\hat H$ on $\FF (H,F)$. Notice that
$\rho (a\circ b)$ acts as $\half E(a\otimes b)$.

Lemma \ref{lemma:bracket} gives $\FF (H,F)$ the structure of a representation of
$\spfrak (H)$. Each graded summand  $\FF (H,F)_i$ is invariant. The 
lemma does not generalize very well to the case when $H$ is infinite dimensional (a case we consider below).  In this respect, the following modification of $\tau$ does a better job.
Suppose we are given an isotropic subspace  ${F'}\subset H$ complementary to $F$.
Then the inclusion ${F'}\subset \hat H$ is a Lie homomorphism and the resulting
map of ${F'}$-modules
\[
\sym_\pt {F'}=U{F'}\to \FF(H,F)
\]
is an isomorphism. In $H\otimes H$ we have a projection operator that
is the identity on $(F\otimes F)\oplus ({F'}\otimes F)\oplus ({F'}\otimes {F'})$ and transposition
on $F\otimes {F'}$. We denote that projector by the \emph{normal ordering} symbol: 
\[
\alpha\in H\otimes H\mapsto\;  :\alpha:\; \in (F\otimes F)\oplus ({F'}\otimes F)\oplus ({F'}\otimes {F'})\subset H\otimes H.
\]
Observe that this restricts to an isomorphism of $\sym^2(H)$ onto
$\sym^2F\oplus ({F'}\otimes F)\oplus \sym^2{F'}$. Although the normal ordering notation is 
in accordance with common practice, it hides the fact that it depends on a choice of 
${F'}$. 

We put
\[
\hat\tau_{F'} : \spfrak(H)\to (U\hat H)_0 [\hbar^{-1}], \quad A\mapsto\hbar^{-1}\widehat{:E^{-1}(A):}
\]
So if $(e_1,\dots ,e_g)$ is a basis of $F$ and $(e_{-1},\dots ,e_{-g})$ is the basis of ${F'}$ 
characterized by $(e_i,e_j)=i\delta_{i+j,0}$, then for $D\in\spfrak(H)$, 
\[
\hat\tau_{F'} (D)=\frac{1}{2\hbar} \sum_{i,j} \frac{(D(e_i),e_j)}{ij}:e_{-i}\circ e_{-j}:\quad , 
\]
where the sum is over all $i,j\in\{-g,\dots ,g\}$.

\begin{lemma}\label{lemma:taudeviation}
If $A\in \spfrak (H)$, then $\hat\tau_{F'}(A)=\tau (A)-\half \trace (A^{{F'}})$, where
$A^{{F'}}\in \End ({F'})$ denotes the restriction of $A$ to ${F'}$
followed by projection along $F$ onto ${F'}$. If moreover $A(F)\subset F$, then left multiplication by $\hat\tau_{F'} (A)$ acts on $\FF(H,F)$ as $A$ on $\sym_\pt F'$.
Finally, $\hat\tau_{F'} :\spfrak (H)\to (U\hat H)_0[\hbar^{-1}]$ is a Lie homomorphism.
\end{lemma}
\begin{proof}
For the first assertion it is enough to verify  this for the  case that $A=E(u\otimes u)$ with $u\in H$. If we  write $u=a+a'$ with $a\in F$, $a'\in {F'}$, then 
\[
:u\circ u:= a\circ a +2a'\circ a +a'\circ a'=
u\circ u- (a\circ a'-a'\circ a)=u\circ u-(a,a')\hbar.
\]
Since $E(u\otimes u)(x)= 2(a,x)a'$, its trace is  $2(a,a')$ and so 
$\hat\tau_{F'}(A)=\tau (A)-\half \trace (A^{{F'}})$. 
If in addition $A$ respects $F$, then $\hat\tau_{F'}(A)\in H\circ F$ and so $\hat\tau_{F'}(A)(v_o)=0$.
If we combine this with the last assertion of Lemma \ref{lemma:bracket} we find that $\hat\tau_{F'} (A)$ acts
on $\FF(H,F)$ as $A$ on $\sym_\pt F'$.

Since $\tau$ is a Lie homomorphism, it follows from the preceding that  $\hat\tau_{F'}$ is also one. 
\end{proof}

The interest of the preceding lemma lies not only in its second property, but also in the 
fact that for infinite dimensional $H$ it may happen 
that $\hat\tau_{F'} (A)$ is defined, but $\tau (A)$ is not, because $A^{{F'}}$ is not of trace class.
\\

\begin{remark}
Any $A$ in the $\spfrak(H)$-stabilizer of $F$ (denoted here by $\spfrak(H)^F$)
induces a transformation in $H/F$ whose trace equals the trace of $A^{F'}$. Since $A$ has trace zero on $H$, this is also minus the trace of $A$ on $F$. So by Lemma  \ref{lemma:taudeviation}, $\hat\tau_{F'}(A)=\tau (A)+\half \trace (A^{F})$ and hence is independent of ${F'}$.
If  $\det (F)$ is regarded as a character of  $\spfrak(H)^F$ and $L$ is half this character,
or more explicitly, if we are given a one dimensional vector space $L$ on which $A\in \spfrak(H)^F$ acts as multiplication by $\half\trace (A^{F})$ plus an intertwining isomorphism $L\otimes L\cong\det (F)$, then  it is clear that this isomorphism identifies the  representation of $\spfrak(H)^F$ defined by $\tau|\spfrak(H)^F$ with that of 
$(\hat\tau_{F'}|\spfrak(H)^F)\otimes L^{-1}$. 
\end{remark}

\section{The case of a polarized weight one Hodge structure}\label{sect:phs}

We take $k=\CC$ from now on. A natural choice for ${F'}$ arises if we are given a fixed real structure on $H$ (which amounts to giving an anti-involution $a\mapsto \bar a$ on $H$) for which $(\; ,\: )$
is defined over $\RR$: $\overline{(a,b)}=(\bar a,\bar b)$) and $F$ has the property that the hermitian form $\sqrt{-1}( a, \bar a)$ is positive definite on $F$. This means $F$ defines a real Hodge structure of weight $1$ on $H$ with $H^{1,0}=F$ and $H^{0,1}=\bar F$, polarized by   $(\; ,\: )$. For we may then take for ${F'}$ the complex conjugate $\bar F$.
Notice that $\bar F$ is indeed totally isotropic for $(\; ,\: )$ and that $\sqrt{-1}( a, \bar a)$ is then negative definite on $\bar F$. We will write $\hat\tau$ for 
$\hat\tau_{\bar F}:\spfrak (H)\to U\hat H [\hbar^{-1}]$.

We  use the isomorphism
\[
\sym_\pt \bar F\to \FF (H,F)
\]
of $\bar F$-modules  to put an inner product on $\FF (H,F)$. First we do this on
$\bar F$ by taking\footnote{In view of the residue formula a more natural choice is
$\la \bar z,\bar w\ra := (2\pi)^{-1}\sqrt{-1}(w, \bar z)$.}
\[
\la \bar z,\bar w\ra := \sqrt{-1}(w, \bar z),\quad  z,w\in F. 
\]
With this 
convention $\la \;,\;\ra$ is indeed positive definite. We extend this in an obvious
way to $\sym_\pt \bar F$ by
\[
\la \bar z_1 \bar z_2\cdots \bar z_n, \bar w_1 \bar w_2\cdots \bar w_m\ra
=
\begin{cases}
\sum_{\sigma \in\Scal_n} \prod_{i=1}^n\la \bar z_i,\bar w_{\sigma (i)}\ra &\text{ if $m=n$ and}\\
0 &\text{otherwise}
\end{cases}
\]
and transfer this inner product to $\FF (H,F)$ via the isomorphism above.
Since $(\; ,\; )$ identifies the dual of  $\bar F$ with  $F$, we may identify the dual
of $\FF (H,F)$ with $\prod_{i=0}^\infty \sym^i F$. The latter may be thought of as
the space $\CC [F^*]$ of formal expansions of functions on $F^*$. It is best 
to think of these as the formal expansions of theta functions on $F^*$ with respect to 
an arbitrary lattice in $F^*$.

\begin{lemma}\label{lemma:adjoint}
Left multiplication by $a\in H$ in $\FF (H,F)$ has as adjoint relative to $\la\; ,\; \ra$
left multiplication by $\sqrt{-1}\bar a\in H$: for all $\alpha,\beta\in \FF (H,F)$ we have
\[
\la \rho (a)\alpha ,\bar \beta\ra=\la \alpha ,\rho (\sqrt{-1} \bar a)\beta\ra.
\]
\end{lemma}
\begin{proof}
First observe that there is no loss in generality in taking $a\in F$, because  interchanging $\alpha$ and $\beta$ and complex conjugating
then also yields the result for $\bar a$.  So we assume $a\in F$. We may (and will) also assume 
that $\alpha$ resp.\ $\beta$  are representable as
products of elements of $\bar F$, say by $\bar z_1\cdots \bar z_n$ resp.\ 
 $\bar w_1\cdots \bar w_m$. If $m\not= n-1$, then both members of the desired
equation vanish and there is nothing to show. So let us take $m=n-1$.
From the identity
\[
\rho (a)\alpha =\sum_{i=1}^n (a,\bar z_i)\bar z_1\cdots\widehat{\bar z_i}\cdots \bar z_n
=-\sqrt{-1}\sum_{i=1}^n \la \bar z_i,\bar a\ra\bar z_1\cdots\widehat{\bar z_i}\cdots \bar z_n.
\]
we easily see that  $\la \rho (a)\alpha,\bar w_1\cdots \bar w_{n-1}\ra$ equals
\[
-\sqrt{-1}\la \bar z_1\cdots \bar z_n, \bar a\bar w_1\cdots \bar w_{n-1}\ra =
\la \bar z_1\cdots \bar z_n, \sqrt{-1}\bar a\bar w_1\cdots \bar w_{n-1}\ra.
\]
 \end{proof}

Let us extend the anti-involution `bar'   in $H$ to the tensor algebra of $H$ by
\[
\overline{a_1\otimes a_2\otimes\cdots \otimes  a_n}:=
\bar a_n\otimes\cdots \otimes \bar a_2\otimes  \bar a_1
\]

\begin{corollary}\label{cor:adjoint}
The adjoint of the action of $\alpha =a_1\otimes a_2\otimes\cdots \otimes  a_n$
on $\FF (H,F)$ is the action of $(\sqrt{-1})^n\bar\alpha$. In particular, 
the 
if $s\in F\otimes F$, then $s+\bar s$ induces in $\FF(H,F)$ an infinitesimal
unitary transformation.
\end{corollary}

\begin{remark}
We see that a vector $a\in H$ acts unitarily in $\FF(H,F)$ if and only if 
$a+\sqrt{-1}\bar a=0$. This means that $a$ is of the form $(1-\sqrt{-1})u$ with $u$ real.
The set of such vectors is a real subspace of $H$ whose complexification is $H$.
Since it plays here the role of a compact real form, we denote it by $H_c$.
If $a, b\in H_c$ and we write $a=(1-\sqrt{-1})u$ and $b=(1-\sqrt{-1})v$, with $u,v\in H_\RR$, then we see that $(a,b)=-2\sqrt{-1}(u,v)\in\sqrt{-1}\RR\hbar$. So $\hat H_c:= H_c+\sqrt{-1}\RR\hbar$ is a real form of $\hat H$. It follows from the preceding that  $\FF(H,F)$ is a unitary representation of $\hat H_c$. It is also irreducible and according to the Stone-von Neumann theorem  it is essentially the only irreducible unitary representation of $\hat H_c$. This implies that for another choice $F'$ of $F$,  there exists an intertwining  isomorphism $\FF (H,F)\to \FF (H,F')$ of inner product spaces. By Schur's lemma, this isomorphism is unique up to a phase factor. This explains in an a priori fashion why we have the projectively flat connection on the universal family of Fock representations of $H$ that we will encounter below.
\end{remark}

The elements of  $\rho (\sym^2 F)$ resp.\ $\rho (\sym^2\bar F)$ clearly 
mutually commute, but an element of one will not in general commute
with an element of the other. In the following section we 
need an expression for the brackets between two such elements.
The isomorphism $E:H\otimes H \to \Hom (H,H)$,  $E(a\otimes b)(x)= 2a (b,x)$, 
has the evident restrictions 
\[
E_F: F\otimes F\to \Hom (\bar F, F), \quad 
E_{\bar F}: \bar F\otimes \bar F\to \Hom (F,\bar F).
\]
Notice that  if $\alpha\in F\otimes F$, then $\overline{E_F(\alpha)}=
E_{\bar F}(\bar \alpha)$.

\begin{lemma}\label{lemma:Tbracket}
If $\alpha, \beta \in \sym^2(F)$, then we have in  $\FF (H,F)$:
\[
[\rho (\bar\alpha) ,\rho (\beta)]=E_{\bar F}(\bar\alpha)E_F(\beta)+\half\trace_{\bar F}(E_{\bar F}(\bar\alpha)E_F(\beta)),
\]
where $E_{\bar F}(\bar\alpha)E_F(\beta)\in \End (\bar F)$ acts in the obvious way
on $\FF (H,F)\cong \sym_\pt(\bar F)$. 
\end{lemma}
\begin{proof}
Since all terms are linear in $\alpha$ and anti-linear in $\beta$ it suffices to check
this for the case $\alpha=a\otimes a$ and $\beta=b\otimes b$ with $a,b\in F$. Recall that we have in $U\hat H[\hbar^{-1}]$
\begin{multline*}
[\bar a\circ \bar a, b\circ b]
=2(\bar a,b)\big(b\circ \bar a+\bar a\circ b\big)\circ\hbar=\\
=4(\bar a,b) \bar a\circ b\circ\hbar+2(\bar a, b)(b,\bar a)\circ\hbar^2.
\end{multline*}
It remains to observe that if $x\in F$, then
\[
E_{\bar F}(\bar a\otimes \bar a)E_F(b\otimes b)(\bar x)=4\bar a (\bar a,b)(b,\bar x)
\]
and that the trace (on $\bar F$) of this transformation equals $4(\bar a, b)(b,\bar a)$.
\end{proof}

\section{Unitary connections on Fock bundles}
Let be given a polarized variation of $\RR$-Hodge structure $(\Hcal,\nabla^\Hcal, \bar{},
(\; ,\; ),  \Fcal)$ of weight one over a complex manifold $S$. We recall that this consists of giving a complex vector bundle $\Hcal$ with a flat connection, an anti-involution
$\,\bar{}\,$ on $\Hcal$ that is flat for $\nabla^\Hcal$ (so that $\Hcal$ is the complexification of 
flat real vector bundle $\Hcal_\RR$),  a flat symplectic form $(\; ,\; )$ on $\Hcal$ that is real relative to $\bar{}$ and  a holomorphic subbundle $\Fcal$ of $\Hcal$  (relative to the holomorphic structure on $\Hcal$ that makes the flat local sections holomorphic). We require that these data define in every fiber  the structure described as in Section \ref{sect:phs}. We then have defined a complex vector bundle $\FF (\Hcal,\Fcal)$ over $S$ of Fock spaces. This bundle is identified with the symmetric algebra  $\sym_\pt(\bar\Fcal)$ and comes with an inner product. We exhibit a natural unitary connection on this bundle.

The connection $\nabla^\Hcal$ defines on the  subbundle $\Fcal$ a \emph{second fundamental form}: if $a$ is a local holomorphic section of $\Fcal$, then consider the projection $\sigma (a)$ of the $\Hcal$-valued $1$-form $\nabla^\Hcal (a)$ in $\bar\Fcal$ along $\Fcal$. This is an element of  $\Ecal^{1,0}(\bar\Fcal )$. If $f$ is a local holomorphic function, then $\sigma (f a)=f\sigma (a)$ and so $\sigma$ defines a bundle map: it is a global section  of $\Ecal^{1,0}(\Hom (\Fcal,\bar\Fcal))$. This section is symmetric relative to
the symplectic form $(\sigma (a),b)=(\sigma (b),a)$. So the identification
\[
E_{\overline{\Fcal}}: \bar\Fcal\otimes\bar\Fcal\to \Hom (\Fcal,\bar\Fcal), \quad
E_{\overline{\Fcal}}(\bar a\otimes \bar b)(x):= 2\bar a (\bar b,x),
\]
makes $\sigma$ correspond to a $1$-form with values in the symmetric part of $\bar\Fcal\otimes\bar\Fcal$. We denote its complex conjugate by $s\in \Ecal^1(\Fcal\otimes\Fcal)$, so that $\bar s=E_{\overline{\Fcal}}^{-1}(\sigma)$. Similarly, we have a  second fundamental form for the subbundle $\bar\Fcal$. Not surprisingly this is $\bar \sigma\in \Ecal^1(\Hom (\bar\Fcal,\Fcal))$, so that $s=E_\Fcal^{-1}(\bar\sigma)$, where
\[
E_{\Fcal}: \Fcal\otimes\Fcal\to \Hom (\bar\Fcal,\Fcal), \quad E_\Fcal(a\otimes b)(x):= 2a (b,x),
\]
So if we are given a local vector field $D$ on $S$, then 
$s(D)$ and $\overline{s(D)}$ are $C^\infty$ sections of $U\hat\Hcal$  and hence
operate in $\FF (\Hcal ,\Fcal)$.

The connection $\nabla^\Hcal$ induces a connections  on the subbundles $\Fcal$ and
$\bar\Fcal$: we let $\nabla^\Fcal$ be the connection on $\Fcal$ obtained
by applying to a section of $\Fcal$ first $\nabla^\Hcal$ and project onto $\Fcal$ along
$\bar\Fcal$. This connection is unitary relative to the restriction of the inner product to $\Fcal$.

The definition of $\nabla^{\bar\Fcal}$ is similar. It is also unitary and is in fact 
equal to the  connection that we get from $(\Fcal, \nabla^\Fcal)$ under the complex 
conjugation map. The connection $\bar\Fcal$ extends in an obvious way to a connection of $\sym_\pt(\bar\Fcal)$. We also view this as a connection on $\FF (\Hcal ,\Fcal)$.

\begin{theorem}\label{thm:Fconnection}
The connection $\nabla^\FF$ in $\FF (\Hcal ,\Fcal)$ defined by 
$\nabla^{\bar\Fcal}+\rho(s +\bar s)$ is unitary and has curvature equal to the scalar
$-\half\trace (\sigma\wedge\bar\sigma)$ (which is also half the 
curvature the connection  $\nabla^{\Fcal}$ induces on the determinant line bundle $\det(\Fcal)$). In particular, the connection
$\nabla^\FF$ is projectively flat.
\end{theorem}
\begin{proof}
Since $\nabla^{\bar\Fcal}$ is unitary in $\bar\Fcal$, it is also unitary in
$\sym_\pt\bar\Fcal$. According to Corollary \ref{cor:adjoint}, 
$\rho (s +\bar s)$ is an infinitesimal unitary transformation of
$\sym_\pt\bar\Fcal$ and it then follows that $\nabla^\FF$ is unitary.

It remains to prove that $\nabla^\FF$ is projectively flat.
The discussion will be local, of course. Choose $o\in S$ and denote the fiber
of $\Fcal$ over $o$ by $F$.  Choose a local $C^\infty$-trivialization of $\Fcal$
over a neighborhood $U$ of $o$ as a \emph{unitary} bundle: So $\Hcal |_U$
has been identified with the constant vector bundle $F_U\oplus \bar F_U$.
Note that this identification takes into account the real structure and the 
symplectic form, but not in general the flat structure. 
The flat connection on $\Hcal_U$ is now given as 
a matrix of $1$-forms, which on $F_U\oplus \bar F_U$ takes the shape
\[
A^\Hcal= \begin{pmatrix}
A^\Fcal & \bar\sigma\\
\sigma & {\overline {A^\Fcal}}\\ 
\end{pmatrix}.
\]
The curvature of $\nabla^\Hcal$ is identically zero and so
\begin{align*}
0&=dA^\Hcal +A^\Hcal\wedge A^\Hcal\\
&=
\begin{pmatrix}
dA^\Fcal +A^\Fcal\wedge A^\Fcal +\bar\sigma\wedge \sigma &
d\bar\sigma + A^\Fcal\wedge\bar\sigma +\bar\sigma\wedge {\overline {A^\Fcal}}\\
d\sigma + \sigma\wedge A^\Fcal +{\overline {A^\Fcal}}\wedge\sigma &
d{\overline {A^\Fcal}} +{\overline {A^\Fcal}}\wedge {\overline {A^\Fcal}} +\sigma\wedge \bar\sigma\\
\end{pmatrix}
\end{align*}
(where the wedge sign between matrix valued forms is to be understood in the 
obvious manner). This amounts to the two identities:
\[
d{\overline {A^\Fcal}} +{\overline {A^\Fcal}}\wedge {\overline {A^\Fcal}} +\sigma\wedge \bar\sigma=0,\quad
d\bar\sigma + A^\Fcal\wedge\bar\sigma +\bar\sigma\wedge {\overline {A^\Fcal}}=0.
\tag{$\dagger$}
\]
The first of these says that the curvature form 
$\Omega(\nabla^{\overline{\Fcal}})$ of $(\bar\Fcal,\nabla^{\overline{\Fcal}})$ 
equals  $-\sigma\wedge \bar\sigma$. Likewise we find that the 
curvature form of $(\Fcal,\nabla^{\Fcal})$ equals  
$-\bar\sigma\wedge\sigma$ so that the curvature of
$\det(\Fcal,\nabla^{\Fcal})$ is $-\trace(\bar\sigma\wedge\sigma)$.

The connection form for $\nabla^\FF$ is given by
\[
A^\FF:={\overline {A^\Fcal}} + s +\bar s,
\]
where we omitted $\rho$ to simplify notation. The curvature form $\Omega(\nabla^\FF)$ of $\nabla^\FF$ is therefore 
\begin{multline*}
\Omega(\nabla^\FF) = d\overline {A^\FF} +\overline {A^\FF} \wedge \overline {A^\FF}\\
=d{\overline {A^\Fcal}}+ ds +d\bar s+
(\overline {A^\Fcal} +s +\bar s)\wedge(\overline {A^\Fcal} + s +\bar s) \\
=(d\overline {A^\Fcal}+\overline {A^\Fcal}\wedge \overline {A^\Fcal}
 +s\wedge \bar s+\bar s\wedge s)+(s\wedge s+\bar s\wedge \bar s) \\
+(ds +\overline {A^\Fcal}\wedge s + s\wedge \overline {A^\Fcal})
+(d\bar s+\overline{A^\Fcal}\wedge \bar s + \bar s\wedge \overline {A^\Fcal}).
\end{multline*}
It is clear that $s\wedge s=0$ and $\bar s\wedge \bar s=0$. Furthermore, 
we may substitute 
$d\overline{A^\Fcal} +
\overline{A^\Fcal}\wedge \overline{A^\Fcal}=-\sigma \wedge\bar\sigma$
and so the theorem follows from the two lemma's below. 
\end{proof}

\begin{lemma}
We have $ds+\overline{A^\Fcal}\wedge s+s\wedge \overline {A^\Fcal}=0$ and  
$d\bar s+\overline{A^\Fcal}\wedge \bar s + \bar s\wedge \overline {A^\Fcal}=0$.
\end{lemma}
\begin{proof}
Let us write 
\[
A^\Hcal= \tilde A+\tilde\sigma\quad\text{with}\quad \tilde A:=
\begin{pmatrix}
A^\Fcal & 0\\
0 &\overline {A^\Fcal}
\end{pmatrix}
\quad\text{and}\quad
\tilde \sigma:=
\begin{pmatrix}
0 & \bar\sigma\\
\sigma & 0
\end{pmatrix}.
\]
and observe that both $\tilde A$ and $\tilde\sigma$ are $\spfrak (H)$-valued $1$-forms on $S$.
We compute
\begin{multline*}
d\tilde\sigma+\tilde A \wedge \tilde\sigma +\tilde\sigma\wedge\tilde A=\\=
\begin{pmatrix}
0 & d\bar \sigma+A^\Fcal\wedge\bar\sigma +\bar\sigma\wedge\overline {A^\Fcal}\\ 
d\sigma+\overline {A^\Fcal}\wedge\sigma +\sigma \wedge A^\Fcal & 0
\end{pmatrix} 
{\buildrel{(\dagger)}\over =}0,
\end{multline*}
where we used ($\dagger$). 
We now apply $\hat\tau$ to both members and let the resulting expression act in $\FF(\Hcal,\Fcal)$. 
Since the coefficients of $\tilde A$ respect $\Fcal$, it follows from Lemma \ref{lemma:taudeviation}
that $\hat\tau (\tilde A)$ acts on $\FF(\Hcal,\Fcal)$ as $\tilde A$ on $\sym_\pt \bar\Fcal$. The latter action is simply that of $\overline {A^\Fcal}$. On the other hand, $\hat\tau (\tilde\sigma)=
\tau (\tilde\sigma)=s+\bar s$ and so we find:
\[
d(s+\bar s)+\overline {A^\Fcal}\wedge (s+\bar s) +(s+\bar s) \wedge\overline {A^\Fcal}.
\]
By taking the homogeneous parts of degree $-2$ and $2$ we get the desired equalities (for $s$ and $\bar s$
respectively).
\end{proof}

\begin{lemma}
We have the following identity of $\End (\FF(\Hcal,\Fcal))$-valued $2$-forms:
$-\sigma\wedge\bar\sigma+s\wedge \bar s+ \bar s\wedge s =-\half \trace 
(\bar\sigma\wedge\sigma))$.
\end{lemma}
\begin{proof}
If $x^1,\dots , x^{2m}$ are local $C^\infty$-coordinates of $S$ at $o$ and $s=\sum_i s_i dx^i$, then
the coefficient of $dx^i\wedge dx^j$, ($i<j$), in $s\wedge \bar s+ 
\bar s\wedge s$ is 
\[
s_i\bar s_j-s_j\bar s_i+\bar s_is_j-\bar s_js_i=[\bar s_i,s_j]-[\bar s_j,s_i].
\]
From Lemma \ref{lemma:Tbracket} we  see that this is equal to
\begin{multline*}
E_{\bar\Fcal}(\bar s_i)E_{\Fcal}(s_j)-E_{\bar\Fcal}(\bar s_j)E_{\Fcal}(s_i)+
\half \trace_{\bar\Fcal}\big(E_{\bar\Fcal}(\bar s_i)E_{\Fcal}(s_j)-(E_{\bar\Fcal}(\bar s_j)E_{\Fcal}(s_i) \big)\\
=(\sigma_i\bar\sigma_j-\sigma_j\bar\sigma_i)
-\half\trace (\bar\sigma_i\sigma_j-\bar\sigma_j\sigma_i),
\end{multline*}
which yields the asserted identity of $2$-forms.
\end{proof}

We can express $\nabla^\FF$ more directly in terms of $\nabla^\Hcal$ as follows.
Let $(e_1,\dots ,e_g)$ be a local basis of holomorphic sections of $\Fcal$ and
extend this to symplectic basis $(e_1,\dots ,e_g,e_{-1},\dots,e_{-g})$ of local sections of 
$\Hcal$: $(e_i,e_j)=\sign (i)\delta_{i+j,0}$. We put
\[
u:=\half \sum_{i,j=1}^g(\nabla^\Hcal (e_j),e_i)e_{-i}\otimes e_{-j}.
\]
Since $(\; ,\; )$ is flat for $\nabla^\Hcal$,  $u$ is symmetric: it is a holomorphic
differential with values in  $\sym^2\Hcal$. We claim that $E_\Fcal(u)=\sigma$. For if $k,l>0$, then
\begin{multline*}
(E_\Fcal(u)(e_k), e_l)=2.\half \sum_{i,j=1}^g(\nabla^\Hcal (e_j),e_i)(e_{-i},e_l) (e_{-j},e_k)\\
=(\nabla^\Hcal (e_k),e_l)=(\sigma (e_k),e_l).
\end{multline*}
This implies that $u$ and $\bar s$ project to the same $\sym^2(\Hcal/\Fcal)$-valued differential on $S$. Hence their images in $U\hat\Hcal$ differ by a $\hbar$-valued differential on $S$. We have equality if  we take $e_{-i}=\sqrt{-1}\bar e_i$.

\begin{proposition}
We have $\nabla^\FF=\nabla^{\Hcal}+\rho(\bar s)$ in the sense that if 
$\alpha=a_1\circ\cdots \circ a_k$, with  $a_1,\dots ,a_k\in\Fcal$, then 
\[\nabla^\FF (\bar \alpha\circ v_o)
=\big(\sum_{i=1}^k \bar a_k\circ\cdots\circ \nabla^\Hcal(\bar a_i)\circ\cdots \circ \bar a_1\circ v_o)\big)+\bar\alpha\circ \bar s\circ v_o
\]
and if we replace $\bar s$ by $u$, then the two members differ by a differential on $S$. 
\end{proposition}
\begin{proof}
Let  $a\in\Fcal$. By Lemma \ref{lemma:bracket} we have 
$[s,\bar a]=\bar\sigma (\bar a)$ in $U\hat\Hcal$. It follows that in $\FF(\Hcal,\Fcal)$ we have the identity
\[
[\nabla^{\overline{\Fcal}}+s, \bar a]=\nabla^{\overline{\Fcal}}(\bar a)+\bar\sigma (\bar a)=\nabla^{\Hcal}(\bar a)
\]
and hence
\begin{multline*}
\nabla^\FF (\bar a_k\circ\cdots \circ \bar a_1\circ v_o)=\\
=\big(\sum_{i=1}^k \bar a_k\circ\cdots\circ \nabla^{\overline{\Fcal}} (\bar a_i)\circ\cdots \circ \bar a_1\circ v_o)\big)+
(s+\bar s)\bar a_k\circ\cdots \circ \bar a_1\circ v_o\\
=\big(\sum_{i=1}^k \bar a_k\circ\cdots\circ \nabla^{\Hcal} (\bar a_i)\circ\cdots \circ \bar a_1\circ v_o)\big)+\bar a_k\circ\cdots \circ \bar a_1\circ \bar s\circ v_o,
\end{multline*}
where we used that  $\bar s$ commutes with $\bar a_i$. The last assertion is clear.
\end{proof}

\begin{remark}
The connection $\nabla^\FF$ becomes flat when twisted with a unitary line bundle:
Suppose we are given  a square root  $\Lcal$ of $\det(\Fcal)$, that is, 
a holomorphic line bundle $\Lcal$ over $S$ equipped with an isomorphism 
$\Lcal^{\otimes 2}\cong \det(\Fcal)$. If we give $\Lcal$ the unitary connection 
$\nabla^\Lcal$ that makes $\nabla^{\Lcal^{\otimes 2}}$
correspond to $\nabla^{\det(\Fcal)}$, then clearly 
$\Omega (\nabla^\Lcal)=\half\trace(\bar\sigma\wedge\sigma)$. 
It follows then from Theorem
\ref{thm:Fconnection} that the tensor product of unitary bundles 
$\FF (\Hcal ,\Fcal)\otimes_{\Ocal_S}\Lcal$ is flat. 

We can also state this differently: if
$p:L\to S$ is the geometric realization of $\Lcal$, then $p^*\Lcal$ has a `tautological'
section whose divisor is the zero section on $p$. So if $p^\circ: L^\circ\to S$ is the restriction to the complement of the zero section (a $\CC^\times$-bundle), then we have 
an identification of  the pull-backs of $\FF (\Hcal ,\Fcal)$ and 
$\FF (\Hcal ,\Fcal)\otimes_{\Ocal_S}\Lcal$ along $p^\circ$ as holomorphic vector bundles. So the resulting (scalar) modification of the connection $\nabla^\FF$ on $(p^\circ)^*\FF (\Hcal ,\Fcal)$ is flat. 
Yet another way to express this fact is to the say that $\FF (\Hcal ,\Fcal)$ becomes a module
over the $\Ocal_S$-algebra of differential operators from $\Lcal$ to itself.
\end{remark}

The preceding suggests the following universal set-up. Let $g:=\half\dim H$ and consider the subspace of the Grassmannian of $g$-dimensional subspaces $F\subset H$ that define  a Hodge structure polarized by $(\; ,\; )$. It has two connected components that are interchanged by complex conjugation. Let $\Scal$ be one of these components. It is an orbit of the real symplectic group 
$\Sp (H_\RR)$  and as such a symmetric space for $\Sp (H_\RR)$ (a stabilizer is a maximal compact subgroup). It is well-known that $\Scal$ is a contractible Stein manifold. Over $\Scal$ we have  a tautological $g$-plane bundle $\Fcal$ (the fiber over $[F]$ is $F$). It is trivial as a holomorphic vector bundle and so the line bundle $\det (\Fcal)$ admits a square root, that is, there exists a holomorphic line bundle $\Lcal$ on $\Scal$ and an isomorphism
$\Lcal^{\otimes 2}\cong\det (\Fcal)$. The choice of such a square root is unique up to isomorphism, the isomorphism being unique up to sign. It defines the \emph{metaplectic group} $\Mp(\Lcal)$: these are the bundle automorphisms of $\Lcal$ that induce in 
$\Lcal^{\otimes 2}\cong\det (\Fcal)$ a transformation in $\Sp (H_\RR)$. It is a connected double cover of $\Sp (H_\RR)$. We denote the central element in the kernel
of $\Mp(\Lcal)\to \Sp (H_\RR)$ by $-1$.

\begin{proposition}
The action of $\spfrak (H_\RR)$ on $\FF (\Hcal,\Fcal)\otimes\Lcal$ defined by 
$\hat\tau$ can be integrated to one of  $\Mp(\Lcal)$ with $-1\in\Mp(\Lcal)$ defining the grading. This action preserves parity.
\end{proposition}

The representation of the metaplectic group on the odd summand $\FF (H,F)_1$ is faithful and is called the \emph{Shale-Weil} representation \cite{lionvergne}.

\begin{remark}
We could introduce in the definition of $\FF (\Hcal ,\Fcal)$ a positive real parameter $\ell$ (a `level') by letting be $\FF_\ell(\Fcal)$ be the representation of $U\hat H$ induced
by the linear map $\hat F=F\oplus\CC\hbar\to\CC$ given by $\ell$ times the coefficient 
of $\hbar$. This however does not really bring us something more general since we may
also obtain this by  replacing $(\; ,\; )$ by $\ell (\; ,\; )$.
As mentioned before, we may think of $\FF (\Hcal ,\Fcal)$ as dual to the space of theta functions on the fibers of $\HH$ realive to some unspecified lattice.  Mumford had observed that such spaces indeed come
with a natural projectively flat connection.
\end{remark}

\section{The WZW connection on a Fock bundle}\label{sect:wzw}
In this section we fix a $\CC$-algebra $R$ and an $R$-algebra $\Ocal$ 
isomorphic to the formal power series ring $R[[t]]$.
In other words, $\Ocal$  comes with a principal ideal $\mfrak$
so that $\Ocal$ is complete for the $\mfrak$-adic topology and $\Ocal/\mfrak^j$ is
for $j=1,2,\dots$ a free $R$-module of rank $j$. We denote by
$K$ the localization of $\Ocal$ obtained by inverting a generator of $\mfrak$.
We denote by $\ord :K\to \ZZ$ the valuation defined by $\mfrak$.
It is clear that for  $N\in\ZZ$, $\mfrak^N\subset K$ is the $\Ocal$-submodule of $K$ of elements of order $\ge N$. The $K$-module of $\CC$-derivations that are continuous (for the $\mfrak$-adic 
topology) from $K$ into $K$ is denoted by $\theta_{K}$ (short for $\theta_{K/\CC}$ since $\CC$ is our base field). The submodule of those that preserve $R$ resp.\ kill $R$ is denoted by $\theta_{K,R}$ resp.\  $\theta_{K/R}$.
The quotient module  $\theta_{K,R}/\theta_{K/R}$ can be identified with 
$\theta_{R}$. If $\omega_{K/R}$ denotes the $K$-dual of $\theta_{K/R}$, then
we have a universal $R$-derivation $d:  K\to \omega_{K/R}$.

These $K$-modules come with filtrations: $F_N\theta_{K,R}$ consists of the derivations that take $\mfrak$ to $\mfrak^{N+1}$ and $F^N\omega_{K/R}$ consists of the $K$-homomorphisms $\theta_{K/R}\to K$ that take  $F^0\theta$ to $\mfrak^N$. In terms of the generator $t$ above, $K=R((t))$, $\theta_{K,R}=\theta_{R}+R((t))\frac{d}{dt}$, $\theta_{K/R}=R((t))\frac{d}{dt}$, $\omega_{K/R} =R((t))dt$ and $F^N\omega_{K/R}=R[[t]]t^{N-1}dt$. 

The residue map
$\res: \omega_{K/R}\to R$ which, in terms of a local parameter $t$, assigns to an element of $R((t))dt$ the coefficient of $t^{-1}dt$
is independent of the choice of $t$. The $R$-bilinear map
\[
 \omega_{K/R}\times K \to R,\quad  (\alpha, f)\mapsto \res (f\alpha)
\]
is a topologically perfect pairing of filtered $R$-modules: $\res (t^k. t^{-l-1}dt)=\delta_{k,l}$ so that any $R$-linear  $\phi: K\to R$  which is continuous (i.e., $\phi$ zero on $\mfrak^N$ for some $N$) is definable by an element of $\omega$
(namely, $\sum_{k<N} \phi (t^k)t^{-k-1}dt$) and likewise
for a $R$-linear continuous map $\omega\to R$. We define an antisymmetric $R$-bilinear form on $K$ by
\[
(f,g):=\res (gdf).
\]
Its kernel is the base ring $R\subset K$, and so we have an induced $R$-symplectic form on $K/R$ (which $d$ maps isomorphically onto $\res^{-1}(0)\subset \omega_{K/R}$). 

We regard $K$ has an abelian $R$-Lie algebra. Its enveloping algebra $UK$
is  the symmetric algebra $\sym^R_\pt K$. We given it  the $\mfrak$-adic topology
and denote by $\Ubar K:=\varprojlim_N UK/\mfrak^N UK$ its completion.

The proof of the following lemma is straightforward.

\begin{lemma} 
Any $D\in\theta_{K/R}$ is,  when viewed as an $R$-linear map $\omega\to K$,
selfadjoint relative to the residue pairing: $\res (\la D,\alpha \ra\beta )
=\res(\la D ,\beta\ra\alpha)$.  
\end{lemma}

Since $D\in\theta_{K/R}$ is selfadjoint we may identify it with an
element of the closure of $\sym^R_2K$ in $\Ubar K$. Let 
$\tau (D)$ be half this element. We shall need  to express this in terms of 
a topological $R$-basis  of $K$. Let us say that a topological $R$-basis of $K$ indexed by $\ZZ$, $(e_i)_{i\in\ZZ}$, is \emph{quasi-symplectic} if
\begin{enumerate}
\item[(i)] $e_0=1\in K$, $e_i\in\mfrak$ for $i>0$, and
\item[(ii)] $(e_i ,e_j)=i\delta_{i+j,0}$ for all $i,j$.
\end{enumerate}
So if $t$ is a generator of $\mfrak$, then $(e_i:=t^i)_{i\in\ZZ}$ is a quasi-symplectic basis
of $K$. If $k$ is a positive integer, then there exists an integer $N_k$ such that
the submodule  $(\sum_{j=1}^{N_k} Re_{-j})\cap \mfrak^{-k}$  maps onto $\mfrak^{-k}/\Ocal$. For $i>N_k$, $e_i$ is perpendicular to this submodule and so must lie in $\mfrak^{k+1}$. This proves that  a quasi-symplectic basis encodes the topology of $K$, because a neighborhood basis of zero  is the collection of  subspaces $\{\prod_{i\ge k} Re_i\}_k$.
In terms of a quasi-symplectic basis, $\tau(D)$ takes the following form
\[
\tau (D)=\frac{1}{2\hbar} \sum_{i,j\in\ZZ-\{ 0\}} \frac{(D(e_i),e_j)}{ij}
e_{-i}\circ e_{-j}.
\]
Observe that the map  $\tau: \theta\to\Ubar K$ is continuous.

The residue map defines a Heisenberg algebra $\hat K$ (also called the
\emph{oscillator algebra}) with underlying  $\CC$-vector space $K\oplus \hbar R$ and Lie bracket 
\[
[f +\hbar r ,g+\hbar s ]:=\res (g\, df)\hbar.
\]
So $[e_k,e_l]=\hbar k\delta_{k,l}$ and the center of  $\hat K$ is 
$R[\hbar,e]$, where $e=t^0$ denotes the unit element of $K$.  It follows that $U\hat K$ is as an $R[\hbar]$-algebra the quotient of the tensor algebra of  $K$ (over $R$) tensored with $R[\hbar]$ by the two-sided ideal generated by the elements $f\otimes g-g\otimes f-(f,g)\hbar$.  The obvious surjection $U\hat K\to UK=\sym_{R,\pt} K$ is the reduction modulo $\hbar$.
 
We complete  $U\hat K$ $\mfrak$-adically on the right:
\[
U\hat K\to \Ubar\hat K:=\varprojlim_N U\hat K /U\hat K\circ \mfrak^N.
\]
Since  for $N\ge 0$, $\mfrak^N$ is an abelian subalgebra of $\hat H$,
this is a completion of $R$-algebra's. Notice that this completion has the collection 
$e_{k_1}\circ\cdots\circ e_{k_r}$ with $r\ge 0$, $k_1\le k_2\le \cdots \le k_r$,
as topological $R[\hbar]$-basis.  (Since $\hat K$ is not abelian, the left and  right
$\mfrak$-adic topologies differ. For instance, 
$\sum_{k\ge 1} e_k\circ e_{-k}$ does not converge in $\Ubar\hat K$, whereas 
$\sum_{k\ge 1} e_{-k}\circ e_k$ does.)
This completion does not affect the Fock representation attached to the pair
$(K,\Ocal)$: $\FF (K,\Ocal)$ is naturally a left $\Ubar\hat K$-module.
\\

A choice of quasi-symplectic basis defines a  (noncanonical) $R$-linear map  
\[
\hat\tau (D)=\frac{1}{2\hbar} \sum_{i,j\in\ZZ-\{ 0\}} \frac{(D(e_i),e_j)}{ij} 
:e_{-i}\circ e_{-j}:\;  \in \Ubar\hat K[\hbar^{-1}],
\]
where  we adhered to the normal ordering convention: the factors inside the
two colons are ordered in nondecreasing order (the factor with highest index 
come last, hence acts first for a left action). So if
$D$ has positive order, then $\hat\tau(D)$ annihilates the generator $v_o$ of 
$\FF(K,\Ocal)$.

\begin{lemma}\label{lemma:chat}
We have $[\hat\tau(D),f]= D(f)$ as an identity  in $\Ubar \hat K[\hbar^{-1}]$ 
(where $f\in K\subset\hat K$) and
 if for $k\not= 0$, $D_k\in\theta_{K/S}$ is characterized by $D_k (e_i)=ie_{k+i}$,
then 
\[
[\hat\tau(D_k), \hat\tau(D_l)]=(l-k)  
\hat\tau(D_{k+l})+\frac{1}{12}(k^3-k)\delta_{k+l,0}.
\]
\end{lemma}
\begin{proof}
For the first statement we compute $[\hat\tau(D_k), e_l]$. The only terms
in the expansion of $\sum_{i+j=k} : e_i\circ e_j:$ that can contribute 
are of the form $[e_{k+l}\circ e_{-l}, e_l]$ or $[e_{-l}\circ e_{k+l}, e_l]$,
(depending on whether $k+2l\le 0$ or $k+2l\ge 0$) and with coefficient 
$\half$ if $k+2l= 0$ and $1$ otherwise. In all cases the result is 
$ le_{k+l}=  D_k(e_l)$.
 
Formula (i) implies that 
\begin{multline*}
[\hat\tau(D_k), \hat\tau(D_l)]=
\lim_{N\to\infty}\sum_{|i|\le N} \frac{1}{2\hbar^2}\left( D_k(e_i)\circ e_{l-i}+ e_i\circ D_k(e_{l-i})\right)\\
=-\hbar^{-1} \lim_{N\to\infty}\sum_{|i|\le N} 
\left( ie_{k+i}\circ e_{l-i}+ e_i\circ (l-i)e_{k+l-i}\right)
\end{multline*}
This is up to a reordering equal to $(l-k)\hat\tau(D_{k+l})$.
The terms which do not commute and are in the wrong order are those
for which $0<k+i=-(l-i)$ (with coefficient $i$) and for which
$0<i=-(k+l-i)$ (with coefficient $(l-i)$). This accounts for
the extra term $\frac{1}{12}(k^3-k)\delta_{k+l,0}$.
\end{proof}

This lemma shows that the transformations $\hat\tau(D)$ generate a $R$-Lie subalgebra 
of $\Ubar \hat K[\hbar^{-1}]$ that is a central (Virasoro) extension if $\theta_{K/R}$ by $R$. We denote this central extension by $\hat\theta_{K/R}$. It is independent of the choice of the quasi-symplectic basis (but the section $\hat\tau :\theta\to \hat\theta$ does depend on it).

It follows from Lemma \ref{lemma:chat} that when $\hat D\in\hat\theta_{K/R}$ 
maps to $D\in\theta_{K/R}$, then
\begin{multline*}
\tag{$\ddagger$}
\hat D (e_{-k_r}\circ\cdots \circ e_{-k_1} v_0)=\\
=\Big(\sum_{i=1}^r e_{-k_r}\circ\cdots\circ
D(e_{-k_i})\circ\cdots \circ e_{-k_1}\Big) v_0+e_{-k_r}
\circ\cdots \circ e_{-k_1}\hat D v_0.
\end{multline*}

We saw that if $\hat D\in\hat\theta$ has positive order, then $\hat Dv_0=0$ and hence
such an element  acts on $\FF(K,\Ocal)$ via $D$ by coefficientwise derivation.

This observation has an interesting consequence.

\begin{corollary}\label{cor:derivevir}
The central extension $\hat\theta_{K/R}$ of  $\theta$ by $R$ naturally extends  to an extension of $R$-Lie algebras $\hat\theta_{K,R}$ of $\theta_{K,R} $ by $R$ in such a manner that
\begin{enumerate}
\item[(i)] the Fock representation of $\hat\theta_{K/R}$ on  $\FF(K,\Ocal)$ extends to $\hat\theta_{K,R}$ such  that it preserves any $U\hat K$-submodule of $\FF(K,\Ocal)$ and \item[(ii)]  if $\hat D\in\hat \theta_{K,R}$ lifts $D\in\theta_{K,R}$, then 
$[\hat\tau(\hat D),f]=Df$ for any $f\in K\subset\hat K$. In particular,
the identity ($\ddagger$) remains valid when $\hat D\in \theta_{K,R}$.
\end{enumerate}
\end{corollary}
\begin{proof}
A quasi-symplectic basis $(e_i)_{i\in\ZZ}$ can be used to split $\theta_{K,R}$:
the set of elements of $\theta_{K,R}$ which kill that basis is a $R$-Lie subalgebra
of  $\theta_{K,R}$ which projects isomorphically onto $\theta_{R}$.
Now if $D\in \theta_{K,R}$, write $D=D_\ver+D_\hor$ with $D_\ver\in \theta_{K/R}$ and 
$D_\hor (e_i)=0$ for all $i$ and define an $R$-linear operator
$\hat D$ in $\FF(K,\Ocal)$ as the sum of $\hat D_\ver$ 
and  coefficientwise derivation by $D_\hor$. This map clearly has the
properties mentioned. As to its  dependence on $t$: it is easily checked that
another choice yields a decomposition of the form $D=(D_\hor+D_0) +(D_\ver -D_0)$  with $D_0$ of positive order and in view of the above $\hat D_0$ then acts in $\FF(K,\Ocal)$ via $D_0$ by coefficientwise derivation.
\end{proof}
 
\begin{definition}
Let us say that a $R$-subalgebra $A\subset K$ is of \emph{Fock  type} if it satisfies the following properties:
\begin{enumerate}
\item[($FT_1$)]  $A$ is as a $R$-algebra flat and of finite type,
\item[($FT_2$)] $A\cap \Ocal =R$ and  the $R$-module $K/(A+\Ocal)$ is  
free of finite rank,
\item[($FT_3$)]  $A$ is totally isotropic relative to the residue form.
\item[($FT_4$)]  Any continuous $R$-derivation of $K$ which preserves $A$ maps 
$A^\perp:=\{ f\in K\,:\,  (f,a)=0 \text{  for all  } a\in A\}$ to $A$. 
 \end{enumerate}
\end{definition}

Let $A$ be an $R$-subalgebra $A\subset K$ of Fock  type. 
It is clear that $A^\perp$ is an $A$-submodule of $K$ which contains $A$. The symplectic pairing induces one on $H_A:= A^\perp/A$:
\[
(a_1 , a_2)_A:=\res (\tilde a_2 d\tilde a_1), \quad \text{where  } \tilde a_i\in A^\perp \text{ lifts  } a_i\in H_A.
\]
We denote by $F_A\subset H_A$ the image of $A^\perp\cap\Ocal$ in $H_A$.

\begin{proposition}\label{prop:fockiso}
The $R$-module $H_A$ is free of rank twice that of $K/(A+\Ocal)$.
The symplectic form on $H_A$ is nondegenerate, has $F_A$ as a maximal isotropic subspace and we have a natural isomorphism 
$\FF (H_A,F_A)\cong\FF (K,\Ocal)_A$.
\end{proposition}
\begin{proof}
Let $g$ the rank of $K/(A+\Ocal)$. Choose  a lift $e_{-1},\dots ,e_{-g}\in K^g$ of a basis of  $K/(A+\Ocal)$ with the property that  $A+\sum_{i=1}^g Re_{-i}$ is totally isotropic.
Notice that $\mfrak$ is also isotropic and that $K^-:=A+\sum_{i=1}^g Re_{-i}$ is a supplement of $\mfrak$ in $K$. Any $R$-linear form on $K^-$ is defined by
taking the symplectic product with some element of $\Ocal$. Let for $i=1,\dots ,g$, 
$e_i\in \mfrak$ be such that $e_i\perp A$ and $(e_i,e_j)=i\delta_{i+j,0}$. Then
$A^\perp=A+\sum_{|i|\le g} Re_i$. So $(e_1,\dots ,e_g)$ resp.\ $(e_{-g},\dots ,e_g)$ projects onto an $R$-basis of $F_A$ resp.\ $H_A$. It is now also clear that the symplectic form is nondegenerate and that $F_A$ is a maximal isotropic subspace.

The natural map
\[
\sym_{R,\pt} K^-\to \FF (K,\Ocal)
\]
is an isomorphism of $K^-$-modules. Taking $A$-coinvariants yields
an isomorphism $\sym_{R,\pt} (K^-)/A\cong \FF (K,\Ocal)_A$. 
We have likewise an isomorphism  $\sym_{R,\pt} (K^-)/A\cong \FF (H_A,F_A)$.
The isomorphism $\FF (H_A,F_A)\cong\FF (K,\Ocal)_A$ thus obtained is easily seen
to be independent of the choice of the $e_i$'s.
\end{proof}

We denote by $\theta_{K,A,R}\subset \theta_{K,R}$
the Lie subalgebra of continuous $\CC$-derivations $K\to K$ that preserve both $R$ and $A$ and by $\hat\theta_{K,A,R}$ its preimage in  $\hat\theta_{K,R}$.

\begin{proposition}\label{prop:A}
The Lie algebra  $\hat\theta_{K,A,R}$ preserves 
$A \circ\Ubar\hat K$ an hence induces a transformation  in $\FF (K,\Ocal)_A$. 
This action of $\hat\theta_{K,A,R}$ on $\FF (K,\Ocal)_A$ factors through
an action of a central extension of  $\theta_{R}$ by $R$.
\end{proposition}
\begin{proof}
The first assertion is clear. Now let $D\in \theta_{K,A,R}\cap \theta_{K/R}$. Then  
according to $(A_4)$, $D$ sends $A^\perp$ to $A$. Since $D$ is also symmetric relative to $(\; ,\; )$,  it then follows that
\[
\hat\tau(D)\in A\circ A^\perp +K\circ \mfrak.
\]
This implies that $\hat\tau(D)(v_o)\in A\circ \FF(K,\Ocal)$ so that $\hat\tau(D)$ induces the zero map
in $H_A$. It follows that the preimage of $\theta_{K,A,R}\cap \theta_{K/R}$
in $\hat\theta_{K/R}$ acts through scalars. This yields the last assertion, because hypotheses ($A_1$) and ($A_2$) imply that any
$\CC$-derivation of $R$ lifts to $\theta_{K,A,R}$
\end{proof}

\begin{remark}
The previous discussion generalizes with little modification to the case where
the $R$-algebra $K$ is isomorphic to a finite direct sum of formal Laurent series:
$K=\oplus_{p\in P} K_p$ with $K_p\cong  R((t))$, where $P$ is a nonempty finite index set. We let $\Ocal=\oplus_{p\in P}\Ocal_p$  and $\mfrak=\oplus_{p\in P}\mfrak_p$. 
If $\res: K\times\omega_{K/R}\to R$ denotes the sum of the residue pairings of 
the summands, then $\res$ is still topologically perfect. In this setting, 
the oscillator algebra $\hat K$ is not the direct sum of the $\hat K_p$, but 
the quotient of  $\oplus_p\hat K_p$ that identifies the central generators 
of the summands with a single $\hbar$. We thus get a Fock representation 
$\FF(K,\Ocal)$ of $\hat K$ that ensures that the unit of every summand $\Ocal_p$ acts  the identity; so it is the induced representation of the rank one representation of $F^0\hat K=\Ocal+R\hbar$ in  $Rv_0$. For an $R$-subalgebra $A\subset K$ of Fock type we find a central extension of $\theta_{R}$ acting on $\FF(K,\Ocal)_A$.
\end{remark}

\section{The Fock construction for a single punctured curve}

Let $C$ be a nonsingular connected complex projective curve of genus $g$ and $P\subset C$ a nonempty finite subset. We put $C^\circ:= C-P$. 
Since $C^\circ$ is affine, the cohomology of its De Rham complex
\[
\begin{CD}
0 @>>> \Ocal (C^\circ)  @>{d}>> \omega (C^\circ) @>>> 0,
\end{CD}
\]
is the complex cohomology of $C^\circ$. In particular,
\[
\omega (C^\circ)/d\Ocal (C^\circ) = H^1(C^\circ; \CC).
\]
The latter fits in the complexification of the short exact sequence
\[
0\to H^1(C)\to H^1(C^\circ)\to \tilde H_0(P)\to 0,
\]
where $\tilde H_0(P)$ stands for the reduced homology of $P$:
the formal linear combinations of elements of $P$ with zero sum.
In terms of the  De Rham complex above, the map 
$H^1(C^\circ;\CC)\to\tilde H_0(P;\CC)$ assigns to an element 
of $\omega (C^\circ)$ its residues at the points of $P$ (recall 
that the sum of the residues is always zero and an exact form
has zero residues). The subspace $\omega (C)\subset \omega (C^\circ)$ maps under this map onto $H^{1,0}(C)$. 

We denote by $(\Ocal_P,\mfrak_P)$ the 
completion of the semi-local ring $\oplus_{p\in P}\Ocal_{C,p }$ with respect to the
ideal defining $P$ and by $K_P$ the localization of $\Ocal$ away 
from $\mfrak$. The natural map
$\Ocal (C^\circ)\to K_P$ is an embedding of algebra's; we denote its image
by $A_P$. The elements of $K_P$ whose differential are the restriction
of an element of $\omega (C^\circ)$  form a subspace $B_P\subset K_P$.
Notice that for every $f\in B_P$, $df$ has zero residues.
It then follows from the preceding that $B_P/A_P$ may be identified
with $H^1(C;\CC)$. 

The universal continuous $\CC$-derivation $d: K_P\to \omega_{K_P}$ has kernel the constants 
$\CC^P\subset K_P$ ($\omega_{K_P}$ is a free $K_P$-module of rank one). The residue map 
\[
R_P:=\sum_{p\in P}\res_p:\omega_{K_P}\to \CC
\]
defines a pairing $\omega_{K_P}\times K_P\to\CC$,
$(\alpha, f)\mapsto R_P(f\alpha)$ which  is nondegenerate in the sense that
it identifies $\omega_{K_P}$ with the continuous dual of $\Ocal_P$. So the resulting
antisymmetric form
\[
(f,g)\in K_P\times K_P\mapsto R_P(gdf):=\sum_{p\in P}\res(gdf)\in\CC
\]
has kernel the space of local constants $\CC^P$. 

The following lemma is a special case of a well-known fact. It connects 
the intersection pairing on $H^1(C;\CC)$ with the residue pairing on $ K_p$.

\begin{lemma}\label{lemma:residue}
If $f_1,f_2\in  d^{-1}\omega(\Ccal^\circ)$, then 
$
\int_C c(f_1)\wedge c(f_2)=-2\pi\sqrt{-1}R_P(f_2df_1)$.
\end{lemma}
\begin{proof}
By assumption $df_i\in\omega_{K_P}$ is the restriction of some 
$\omega_i\in\omega(C^\circ)$. For $p\in P$, let  $z_p$ be a local analytic coordinate for $C$ at $p$ whose domains are pairwise disjoint. We assume that $z_p$  maps its domain isomorphically onto the complex unit disk. Let $\phi_p$ be a $C^\infty$ function on $[0,\infty)$ with support inside $[1,0)$ and constant $1$ on $[0,\half]$. Then 
$\tilde\phi_p:=\phi_p (|z_p|)$ is $C^\infty$ on $C$, when extended by zero to all of $C$.
We put $\tilde\phi:=\sum_{p\in P}\tilde\phi_p$.
The form $\tilde\omega_i:=\omega_i-d(\tilde\phi f_i)$ is zero on a neighborhood of $P$
and is closed. Since it has the same periods as $\omega_i$, it represents $c(f_i)$ as a De Rham cohomology class. It follows that
\[
\int_C c(f_1)\wedge c(f_2)=\int_C \tilde\omega_1\wedge \tilde\omega_2
\]
In the right hand side we may replace $\tilde\omega_1$ by $\omega_1$,
because the difference $d(-\tilde\phi f_1)\wedge\tilde\omega_2=
d(-\tilde\phi f_1\tilde \omega_2)$ is exact as a $2$-form on $C$. So the integral becomes
\begin{multline*}
\int_C\omega_1\wedge \tilde\omega_2=\int_{C^\circ} \omega_1\wedge d(-\tilde\phi f_2)=
\int_{C^\circ} d(\tilde\phi f_2\omega_1)=\\
=\lim_{r\downarrow 0}\sum_{p\in P}-\int_{|z_p|= r} \tilde\phi_p f_2\omega_1=
\lim_{r\downarrow 0}\sum_{p\in P}\int_{|z_p|= r}  -f_2\omega_1=\\
=\sum_{p\in P}-2\pi\sqrt{-1}\res_p(f_2\omega_1)=-2\pi\sqrt{-1}R_P(f_2df_1).
\end{multline*}
\end{proof}

\begin{corollary}
The subalgebra $A_P$ of $K_P$ is of Fock type with $A_P^\perp=B_P$ and $A_P^\perp/A_P$ is canonically identified $H^1(C;\CC)$ in such a manner that the intersection pairing on the latter is represented by
$-2\pi\sqrt{-1}$ times the residue pairing. Moreover $F^1A_P^\perp=A_P^\perp\cap\mfrak_P$ maps isomorphically onto $H^{1,0}(C)$.
\end{corollary}
\begin{proof}
If $f\in K_P$ is such that $df$ is the image of some $\alpha\in\omega(C^\circ)$ and $g\in A_P$ comes from $\tilde g\in \Ocal(C^\circ)$, then $R_P(gdf)$ is the residue sum  of
an element of $\omega(C^\circ)$ and hence zero. This proves that $B_P\subset A_P^\perp$. 
Lemma \ref{lemma:residue} implies that the residue form induces a nondegenerate pairing on $B_P/A_P$. This implies that $B_P=A_P^\perp$. The other assertions are clear.
\end{proof}

Denote by $ K_P^-\subset  B_P$ the preimage  of 
$H^{0,1}(C;\CC)$ under the surjection $B_P\to H^1(C;\CC)$.  
Then the natural map 
\[
 K_P^-/A_P\to  K_P/(A_P+\mfrak_P)
\] 
is an isomorphism (the former maps isomorphically onto $H^{0,1}(C)$, the latter maps isomorphically onto $H^1(C,\Ocal_C)$ with the map between them being the
natural identification $H^{0,1}(C)\cong H^1(C,\Ocal_C)$). It follows that the natural map  
\[
K^-_P\to  K_P/\mfrak_P
\]
is an isomorphism. Observe that  $K^-_P$ is also totally isotropic for the residue pairing.
These properties make it possible and  tractable to choose the coefficients in the
Fock representation attached to $K_P$  in $K^-_P$. We cannot, however, expect that $K^-_P$ be multiplicatively closed.

\begin{example}
Let $C$ be a hyperelliptic curve of genus $g>0$ and $p\in C$ a Weierstra\ss\ point.
We suppose  $C^\circ=C-\{ p\}$ given as the affine plane curve $y^2=f(x)$ with
$f\in \CC[x]$ a polynomial  of degree $2g+1$ with distinct roots. The rational function $x^g/y$ may serve as a local parameter at $p$, but we prefer to take a $t\in K_p$ for which 
$x=t^{-2}$ and $y=\sqrt{f(t^{-2})}=t^{-1-2g}u(t^2)^{-1}$, 
where $u(t):=(t^{2g+1}f(t^{-1}))^{-1/2}$ is a unit in $\CC[[t]]$. We have
\begin{align*}
A_p&=\Ocal (C^\circ)= \CC[x]+\CC[x]y=\CC[t^{-2}]+\CC[t^{-2}]t^{-1-2g}u(t^2)^{-1},\\
\omega (\Ccal^\circ)&=\CC[x]y^{-1}dx +\CC[x]dx=\CC[t^{-2}]t^{2g-2}u(t^{2})dt+\CC[t^{-2}]t^{-3}dt,\\
\omega (\Ccal)&=\la 1,x,\dots ,x^{g-1}\ra y^{-1}dx=
\la1,t^2, t^4,\dots t^{2g-2}\ra u(t^{2})dt.
\end{align*}
It follows that $B_p$ is spanned by $A_p$ and the $\{\phi_{2i-1}\in K_p\}^{g}_{i=1-g}$, where $\phi_{2i-1}\in t^{2i-1}\CC[[t]]$ is such that $\phi_{2i-1}'=u(t^2)t^{2i-2}$ so that 
$K^-_p$ is the span of $A_p$ and $\phi_{-1}, \phi_{-3},\dots ,\phi_{1-2g}$.
One may check that $K^-_p$ is not multiplicatively closed.
\end{example}

Let us denote by $\hat K_P$ the Heisenberg algebra defined by residue pairing (with underlying vector space $ K_P\oplus\CC\hbar$) and by $\FF ( K_P,\Ocal_P)$ the Fock representation 
(induced by the linear form on $\hat\Ocal_P=\Ocal_P\oplus\CC\hbar$ defined by the coefficient of $\hbar$).
Then the map
\[
\sym_\pt  K^-_P\to \FF ( K_P,\Ocal_P)
\]
is an isomorphism of $K^-_P$-modules. By  taking on both sides the covariants relative to the abelian subalgebra $A_P$ of $\hat K_P$, we get an isomorphism
from $\sym_\pt H^{0,1}(C,\CC)$ onto $\FF (K_P,\Ocal_P)_{A_P}$. The resulting isomorphism 
\[
\FF(H^1(C,\CC),H^0(\omega_S))\cong  \FF (K_P,\Ocal_P)_{A_P}.
\]
can however be obtained without reference to $H^{0,1}(C)$ as  the maps
\[
K_P\supset B_P\to  B_P/A_P{\buildrel d\over \cong} 
H^1(C,\CC), \quad 
\Ocal_P\supset  B_P\cap\mfrak_P {\buildrel d\over \cong} F,
\]
induce maps
\[
 \FF ( K_P,\Ocal_P)_{A_P}\leftarrow
 \FF ( B_P,\Ocal_{\Pcal}\cap B_P)_{A_P}\to \FF(H^1(C,\CC),H^0(\omega_S))
\]
which are both isomorphisms.
\\

Suppose  that $\pi:\tilde C\to C$ is a connected  nonsingular  projective curve over $C$ such that there is no ramification over $P$. If $\tilde P$ denotes the preimage of $P$, then (via $\pi^*$) we have $K_P\subset K_{\tilde P}$, $A_P=A_{\tilde P}\cap K_P$, $B_P=B_{\tilde P}\cap K_P$, $B_P\cap \mfrak_{\tilde P}=
B_P\cap \mfrak_P$ and $R_{\tilde P}| \omega (C^\circ)=\deg (\pi).R_P$.
The obvious map
\[
B_P/A_P\to B_{\tilde P}/A_{\tilde P}
\]
corresponds to the pull-back $\pi^* :H^1(C;\CC)\to  H^1(\tilde C;\CC)$. The preimage
of $H^{0,1}(\tilde C)$ is $H^{0,1}( C)$ and so $K^-_{\tilde P}\cap B_P=K^-_P$.
It follows that we have a natural embedding
\[
\FF (K_P,\Ocal_P)_{A_P}\to \FF (K_{\tilde P},\Ocal_{\tilde P})_{A_{\tilde P}},
\]
which can be understood as the embedding
\[
\sym_\pt H^{0,1}(C)\to \sym_\pt H^{0,1}(\tilde C).
\]
Since this map multiplies the inner product by $\deg (\pi)$, we replace it by its 
multiple by $\deg (\pi)^{-1}$, so that the embedding becomes isometric. If the 
covering is Galois with Galois group $G$, then $G$
acts on $H^{0,1}(\tilde C)$ and  $H^{0,1}(\tilde C)^G$ is the image of  $H^{0,1}(C)$.

For a fixed $C$ and $P$, the  nonsingular connected  covers $\tilde C\to C$ 
that are not ramified over $P$ form a direct system and so we find a Fock space
\[
\tilde\FF (K_P,\Ocal_P):= \varinjlim\FF (K_{\tilde P},\Ocal_{\tilde P}),
\]
and an isometric embedding
\begin{multline*}
\sym_\pt H^{0,1}(C)\cong \FF (K_P,\Ocal_P)_{A_P}\to \\
\varinjlim\FF (K_{\tilde P},\Ocal_{\tilde P})_{\tilde A_P}\cong  \varinjlim\sym_\pt H^{0,1}(\tilde C)\cong \varinjlim\sym_\pt  (K^-_{\tilde P}/A_{\tilde P}).
\end{multline*}
Since any smooth projective curve is a finite cover of $\PP^1$, any two 
connected smooth projective curves are dominated by a connected common finite cover
(take a connected component of the normalization of their fiber product over $\PP^1$).
So this limit is huge (it does not have a countable basis).

\section{The Fock construction for a family of punctured curves}

We here assume given a family $\pi :\Ccal \to S$ of curves ($\pi$ that  is smooth and proper with connected genus $g$ curves as fibers) and a finite nonempty union of
mutually  disjoint sections  $\Pcal\subset\Ccal$. We shall also assume that $S$ is nonsingular and connected  (so $\Ccal$ will then also  have these properties). We put $\Ccal^\circ:=\Ccal-\Pcal$. The previous construction applies fiberwise: the local system
$\HH_\RR:=R^1\pi_*\RR_\Ccal$ comes with a symplectic form. 
Its underlying holomorpic vector bundle 
\[
\Hcal:=\HH_\RR \otimes_{\RR_S}\Ocal_S= R^1\pi_*\pi^{-1}\Ocal_S.
\]
comes with a  flat connection $\nabla^\Hcal$ (whose flat sections define $\HH:=\HH_\RR\otimes\CC\subset\Hcal$) and a  flat symplectic form.
The De Rham complex defines a short exact sequence of holomorphic vector bundles
\[
0\to \pi_*\omega_{\Ccal/S}\to \Hcal\to R^1\pi_*\Ocal_\Ccal\to 0.
\]
We put $\Fcal:= \pi_*\omega_{\Ccal/S}$. Then the $C^\infty$-subbundle $\bar\Fcal$
of $\Hcal$ maps isomorphically onto $R^1\pi_*\Ocal_\Ccal$. As in the absolute case,
we find that $\Acal_P$ is a sheaf of Fock subalgebras of $\Kcal_P$ and that we have
an isomorphism of  vector bundles
\[
\FF(\Hcal,\Fcal)\cong \FF (\Kcal_\Pcal,\Ocal_\Pcal)_{\Acal_\Pcal},
\]
where $\Acal_\Pcal:=\pi_*\Ocal_{\Ccal^\circ}$. It is induced by the homomorphisms of $\Ocal_S$-modules
\[
\Kcal_{\Pcal}\supset\Bcal_{\Pcal}\to \Bcal_{\Pcal}/\Acal_{\Pcal}\cong
\Hcal,\quad 
\Ocal_{\Pcal}\supset \Bcal_{\Pcal}\cap\Ocal_{\Pcal}\to (\Bcal_{\Pcal}\cap\Ocal_{\Pcal})/\Ocal_S \cong \Fcal,
\]
from which we see that this isomorphism is holomorphic. 

We wish to compare the WZW-connection on the  left with the connection defined above on the right.

Denote by $\theta_{\Ccal,S}$ the sheaf of $\Ocal_{\Ccal}$ into itself
that cover a $\pi^{-1}\Ocal_S$-derivation, or what amounts to the sheaf of vector fields
on $\Ccal$ that lift a vector field on $S$. So we have an exact sequence
\[
0\to \theta_{\Ccal/S}\to \theta_{\Ccal,S}\to \pi^{-1}\theta_S\to 0.
\]
Notice that $\pi_*\theta_{\Ccal^\circ,S}$
acts on $\Kcal_{\Pcal}$.

\begin{lemma}\label{lemma:covder}
The action of $\pi_*\theta_{\Ccal^\circ,S}$ on $\Kcal_{\Pcal}$ preserves the 
$\Ocal_S$-submodules $\Acal_{\Pcal}$ and $\Bcal_{\Pcal}$
and the action of $\pi_*\theta_{\Ccal^\circ,S}$ on their quotient $\Bcal_{\Pcal}/\Acal_{\Pcal}$
factors through $\theta_S$. This makes $\Bcal_{\Pcal}/\Acal_{\Pcal}$ a $\Dcal_S$-module which via  the isomorphism $\Bcal_{\Pcal}/\Acal_{\Pcal}\cong\Hcal$
is just covariant derivation.
\end{lemma}
\begin{proof}
Let $\tilde D\in\pi_*\theta_{\Ccal^\circ,S}$. If $f$ is a local section of $\pi_*\Ocal_{\Ccal^\circ}$, then
clearly, so is $\tilde D(f)$. If $f$ is a local section of $\Bcal_{\Pcal}$, then $d_{\Kcal_{\Pcal}/S}f$ is the restriction of a local section $\alpha$ of $\pi_*\omega_{\Ccal^\circ/S}$. Locally on $S$, $\alpha$ lifts to
an absolute differential $\tilde\alpha\in\pi_*\omega_{\Ccal^\circ}$. The Lie derivative $L_{\tilde D}(\tilde\alpha)$ also lies in $\pi_*\omega_{\Ccal^\circ}$. Near $\Pcal$, this Lie derivative equals $d_{\Kcal_{\Pcal}/S}\tilde D(f)$, which shows that $\tilde D(f)\in \Bcal_{\Pcal}$.
 
If $D=0$, then there is no need to lift $\alpha$ to an absolute differential.
Since $\alpha$ is relatively closed, $L_{\tilde D}(\alpha)$ is the relative differential 
of the local section $\la \tilde D, \alpha\ra$ of $\pi_*\omega_{\Ccal^\circ/S}$ and hence
$\tilde D(f)\in\la \tilde D, \alpha\ra$. It follows that $\tilde D(f)\in\Acal_{\Pcal}$. 

A standard argument \cite{deligne} identifies the resulting action of $\Dcal_S$ on $\Bcal_{\Pcal}/\Acal_{\Pcal}$
with covariant derivation in $\Hcal$.
\end{proof}

\begin{theorem}\label{thm:fockiso}
The isomorphism  $\FF(H^1(C,\CC),H^0(\omega_S))\cong  \FF (K_{\Pcal},\Ocal_{\Pcal})_{\pi_*\Ocal_{C^\circ}}$ takes the Fock connection $\nabla^\FF$ to the  WZW-connection
and thus makes the latter a unitary connection.
\end{theorem}
\begin{proof}
As explained in Section \ref{sect:wzw}, we have a central extension 
$\hat\theta_{\Ccal^\circ,S}$ 
of $\theta_{\Ccal^\circ,S}$ by $\Ocal_S$ which acts on $\FF(\Kcal_{\Pcal},\Ocal_{\Pcal})$
according to the following rule: if $\hat D\in \hat\theta_{\Ccal^\circ,S}$ 
lifts $D\in\theta_{\Ccal^\circ,S}$, then
\[
\hat D (f_r\circ\cdots \circ f_1\circ v_o)
=\Big(\sum_{i=1}^r f_r\circ\cdots
D(f_i)\circ\cdots \circ f_1\Big)\circ v_o+f_r\circ\cdots \circ f_1\circ \hat D v_o,
\]
where it remains to explain $\hat D v_o$. We here describe the latter element up
to $\Ocal_S$-multiple of $v_o$ in terms of $D$ only; for a more complete discussion we refer to Section \ref{sect:wzw}. Choose  a local basis
$(\omega_1,\dots ,\omega_g)$ of $\pi_*\omega_{\Ccal/S}$. Let $e_j\in\mfrak_{\Pcal}$ 
be such that $de_j=j\omega_j$. We extend this  to a quasi-symplectic basis
of $K_\Pcal$ in the sense of Section \ref{sect:wzw} that is well adapted to the situation:
Choose $\{e_{-i}\in \Bcal_{\Pcal}\}_{i=1}^g$ such that 
$\res (e_{-i}\omega_j)=j\delta_{ij}$ and $\res (e_{-i}de_{-j})=0$ (here $i,j>0$)
and put $e_0=1\in \Kcal_{\Pcal}$. Let $\{e_{-i}\}_{i=g+1}^\infty$ be a basis of 
$\Acal_{\Pcal}=\pi_*\Ocal_{\Ccal^\circ}$.
Then $\{ e_i\}_{i\le 0}$ is a basis of $\Kcal^-_{\Pcal}$ and $\{ e_i\}_{i\le g}$ is basis of $\Bcal_{\Pcal}$. We extend this to a quasi-symplectic topological basis $\{ e_i\in \Kcal_{\Pcal}\}_{i\in\ZZ}$ of $\Kcal_{\Pcal}$. Now
$\hat D v_o$ has the property that if we put
\[
\hat\tau (D):=-\pi\sqrt{-1}\sum_{i,j}\res \big(\frac{e_jd(De_i)}{ij}\big):e_{-i}\circ e_{-j}: \; \in \Ubar\hat K_{\Pcal},
\] 
then 
\[
\hat D v_o\equiv \hat\tau (D)\circ v_o \pmod{\Ocal_Sv_o}.
\]
It is clear that in $ \hat\tau (D)v_o$ only the terms with $i, j\ge 0$ contribute. Since
$D$ preserves $\pi_*\Ocal_{\Ccal^\circ}$, it induces a transformation in 
$\FF(\Kcal_{\Pcal},\Ocal_{\Pcal})_{\pi_*\Ocal_{\Ccal^\circ}}$. The terms with $i,j>g$ act trivially in   $\FF(\Kcal_{\Pcal},\Ocal_{\Pcal})_{\pi^*\Ocal_{\Ccal^\circ}}$ and since for $i, j=1,\dots ,g$, 
\[
\res (e_jd(De_i))=\res (e_jd\la D, \omega_i\ra)=-\res (\la D, \omega_i\ra \omega_j),
\]
it follows that the image of $\hat\tau (D)v_o$ in $\FF(\Kcal_{\Pcal},\Ocal_{\Pcal})_{\pi_*\Ocal_{\Ccal^\circ}}$
equals
\[
\pi\sqrt{-1}\sum_{i,j=1}^g \res \big(\frac{\la D, \omega_i\ra\omega_j}{ij}\big)e_{-i}\circ e_{-j}\circ v_o.
\]
Up to a scalar in $\Ocal_S$, the last action only depends on $D_S$. 
It follows that under the isomorphism $\FF(\Kcal_{\Pcal},\Ocal_{\Pcal})_{\pi_*\Ocal_{\Ccal^\circ}}\cong \FF(\Hcal,\Fcal)$, this action coincides with the one on $\FF(\Hcal,\Fcal)$. 
\end{proof}

\end{document}